\documentclass[12pt]{amsart}
\usepackage{amscd,amsmath,amsthm,amssymb}
\usepackage{amsfonts,amssymb,amscd,amsmath,enumerate,verbatim}
\usepackage[left]{lineno}
\usepackage{pstricks, pst-plot,pst-3d}
\usepackage{tikz}
\newpsstyle{fatline}{linewidth=1.5pt}
\definecolor{verylight}{gray}{0.97}
\definecolor{light}{gray}{0.9}
\definecolor{medium}{gray}{0.85}
\definecolor{dark}{gray}{0.6}
\def\Gc{{\mathcal G}}

\def\Gc{{\mathcal G}}

\def\opn#1#2{\def#1{\operatorname{#2}}} 
\opn\chara{char} \opn\length{\ell} \opn\pd{pd} \opn\rk{rk}
\opn\projdim{proj\,dim} \opn\injdim{inj\,dim} \opn\rank{rank}
\opn\depth{depth} \opn\grade{grade} \opn\height{height}
\opn\embdim{emb\,dim} \opn\codim{codim}
\opn\Cl{Cl}

\opn\Tr{Tr} \opn\bigrank{big\,rank}
\opn\superheight{superheight}\opn\lcm{lcm}
\opn\trdeg{tr\,deg}
	\opn\reg{reg} \opn\lreg{lreg} \opn\ini{in} \opn\lpd{lpd}
	\opn\size{size} \opn\sdepth{sdepth}
	\opn\link{link}\opn\fdepth{fdepth}\opn\lex{lex}
	\opn\tr{tr}\opn\del{del}
	\opn\type{type}
	\opn\gap{gap}
	\opn\arithdeg{arith-deg}
	\opn\revlex{revlex}
	\opn\div{div} \opn\Div{Div} \opn\cl{cl} \opn\Cl{Cl}
	\opn\Spec{Spec} \opn\Supp{Supp} \opn\supp{supp} \opn\Sing{Sing}
	\opn\Ass{Ass} \opn\Min{Min}\opn\Mon{Mon}
	\opn\Ann{Ann} \opn\Rad{Rad} \opn\Soc{Soc}
	\opn\Im{Im} \opn\Ker{Ker} \opn\Coker{Coker} \opn\Am{Am}
	\opn\Hom{Hom} \opn\Tor{Tor} \opn\Ext{Ext} \opn\End{End}
	\opn\Aut{Aut} \opn\id{id}
	
	\opn\nat{nat}
	\opn\pff{pf}
	\opn\Pf{Pf} \opn\GL{GL} \opn\SL{SL} \opn\mod{mod} \opn\ord{ord}
	\opn\Gin{Gin} \opn\Hilb{Hilb}\opn\sort{sort}
	\opn\PF{PF}\opn\Ap{Ap}
	\opn\mult{mult}
	\opn\bight{bight}
	\opn\div{div}
	\opn\Div{Div}
	\opn\aff{aff}
	\opn\relint{relint} \opn\st{st}
	\opn\lk{lk} \opn\cn{cn} \opn\core{core} \opn\vol{vol}  \opn\inp{inp} \opn\nilpot{nilpot}
	\opn\link{link} \opn\star{star}\opn\lex{lex}\opn\set{set}
	\opn\width{wd}
	\opn\Fr{F}
	\opn\QF{QF}
	\opn\G{G}
	\opn\type{type}\opn\res{res}
	\opn\conv{conv}
	\opn\Deg{Deg}
	\opn\Sym{Sym}
	\opn\Con{Con}
	\opn\gr{gr}
	
	\def\pot#1#2{#1[\kern-0.28ex[#2]\kern-0.28ex]}

	\opn\dirlim{\underrightarrow{\lim}}
	\opn\inivlim{\underleftarrow{\lim}}
	\let\to=\rightarrow
	
	\def\Implies{\ifmmode\Longrightarrow \else
		\unskip${}\Longrightarrow{}$\ignorespaces\fi}
	\def\implies{\ifmmode\Rightarrow \else
		\unskip${}\Rightarrow{}$\ignorespaces\fi}
	\def\iff{\ifmmode\Longleftrightarrow \else
		\unskip${}\Longleftrightarrow{}$\ignorespaces\fi}

	\let\:=\colon
	\newtheorem{Theorem}{Theorem}[section]
	\newtheorem{Lemma}[Theorem]{Lemma}
	\newtheorem{Corollary}[Theorem]{Corollary}
	
	\newtheorem{Remark}[Theorem]{Remark}

	\let\epsilon\varepsilon
	\let\kappa=\varkappa
	\def\qed{\ifhmode\textqed\fi
		\ifmmode\ifinner\quad\qedsymbol\else\dispqed\fi\fi}
	\def\textqed{\unskip\nobreak\penalty50
		\hskip2em\hbox{}\nobreak\hfil\qedsymbol
		\parfillskip=0pt \finalhyphendemerits=0}
	\def\dispqed{\rlap{\qquad\qedsymbol}}
	
	\opn\dis{dis}
	\def\pnt{{\raise0.5mm\hbox{\large\bf.}}}
	
	\opn\Lex{Lex}

\begin{document}
\title[Algebras with straightening laws]{An Algebra with straightening laws which is an integral domain and which is not Cohen--Macaulay}

\author[T.~Hibi]{Takayuki Hibi}
\author[A. A. Qureshi]{Ayesha Asloob Qureshi}
\address{(Takayuki Hibi) Department of Pure and Applied Mathematics, Graduate School of Information Science and Technology, The University of Osaka, Suita, Osaka 565--0871, Japan}
\email{hibi@math.sci.osaka-u.ac.jp}
\address{(Ayesha Asloob Qureshi) Faculty of Engineering and Natural Sciences, Sabanci University, Orta Mahalle, Tuzla 34956, Istanbul, Turkey}	
\email{aqureshi@sabanciuniv.edu, ayesha.asloob@sabanciuniv.edu}



\begin{abstract}
A $4$-dimensional algebra with straightening laws which is an integral domain and which is not Cohen--Macaulay over an arbitrary field is constructed.  Furthermore, it is shown that every three-dimensional graded ASL domain is Cohen--Macaulay. 
\end{abstract}	
\maketitle
\thispagestyle{empty}

\section*{Introduction}
The research area so-called {\em Combinatorics and Commutative Algebra} was born in the middle of 1970s out of Stanley \cite{UBC} and Reisner \cite{Rei} under the great influence of Hochster \cite{Hoc}.  In the early of 1980s, with the appearance of Cohen--Macaulay and shellable partially ordered sets \cite{Bac, Bjo}, Eisenbud \cite{Eis} introduced a new concept of algebras with straightening laws.  We abbreviate an algebra with straightening laws as an ASL.  In the early stage on the study of ASL, it was conjectured that every ASL domain is Cohen--Macaulay.  It is shown \cite[Corollary, p.~32]{HW1} that every three-dimensional homogeneous ASL domain is Cohen--Macaulay.  On the other hand, when the characteristic of the base field is $2$, Terai \cite{Terai} constructs a $4$-dimensional homogeneous ASL domain which is not Cohen--Macaulay.  

In the present paper, to supply the above old results, we show the following

    \begin{Theorem}
\label{ASLdomain}
Every $three$-dimensional graded ASL domain is Cohen--Macaulay.   
\end{Theorem}

\begin{Theorem}
\label{counterexample}
There exists a $4$-dimensional graded ASL domain which is not Cohen--Macaulay over an arbitrary field.   
\end{Theorem}

After recalling fundamental materials on ASL in Section $1$, Theorem \ref{ASLdomain} is proved in Section $2$.  In section $3$, a simple and effective construction of graded ASL is given (Lemma \ref{gradedASL}).  Finally, in Section $4$, to prove Theorem \ref{counterexample}, a $4$-dimensional graded ASL domain which is not Cohen--Macaulay over an arbitrary field is constructed.

We refer the reader to \cite{HIBIred, StaGREEN} for fundamental materials on partially ordered sets and simplicial complexes and to \cite[Chapter 2]{HHgtm260} for the introduction to Gr\"obner bases.

\section{Algebras with straightening laws} 
Every poset (partially ordered set) to be considered is finite.  A {\em chain} of $P$ is a totally ordered subset of $P$.  The {\em length} of a chain $C$ of $P$ is $|C|-1$.  The {\em rank} of $P$ is the maximum length of chains of $P$.  Let $\rank P$ denote the rank of $P$.  We say that a poset $P$ is {\em pure} if every maximal chain of $P$ has the same length.  A pure poset $P$ of rank $d$ is {\em strongly connected} if $C$ and $C'$ are maximal chains of $P$, then there is a sequence $C_0,C_1,\ldots, C_q$ of maximal chains of $P$ with $C_0=C$ and $C_q=C'$ for which $|C_{i-1} \cap C_{i}| = d$ for each $1 \leq i \leq q$.

Let $R$ be a noetherian algebra over a field $K$.  Let $P$ be a poset and suppose that an injection $\varphi: P \hookrightarrow R$
for which the $K$-algebra $R$ is generated by $\varphi(P)$ over $K$ is given.  A {\em standard monomial} is an element of $R$ of the form $\varphi(\gamma_1) \varphi(\gamma_2)\cdots \varphi(\gamma_s)$, where $\gamma_1 \leq \gamma_2 \leq \cdots \leq \gamma_s$ in $P$.  We call $R$ an {\em algebra with straightening laws} \cite{Eis} on $P$ over $K$ if the following conditions are satisfied:
\begin{itemize}
\item[]
(ASL\,-1)
The set of standard monomials is a $K$-basis of $R$;
\item[]
(ASL\,-2)
If $\alpha$ and $\beta$ in $P$ are incomparable and if
\begin{eqnarray*}
\label{ASL}
\, \, \, \, \, \, \, \, \, \, \, \, \, \, \, \, \, \, \, \, 
\varphi(\alpha)\varphi(\beta) 
= \sum_{i} r_i\,\varphi(\gamma_{i_1})\varphi(\gamma_{i_2}) \cdots , \, \, \, 0 \neq r_i \in K, \, \, \, \gamma_{i_1}\leq \gamma_{i_2}\leq \cdots, 
\end{eqnarray*}
is the unique expression for $\varphi(\alpha)\varphi(\beta) \in R$ as a linear combination of distinct standard monomials guaranteed by (ASL\,-1), then $\gamma_{i_1} \leq \alpha, \beta$ for every $i$. 
\end{itemize}
The relations in (ASL\,-2) are called the {\em straightening relations} for $R$.  Note that the right-hand side of the relation in (ASL\,-2) is allowed to be the empty sum $(=0)$, but that, though $1$ is a standard monomial, no $\varphi(\gamma_{i_1})\varphi(\gamma_{i_2}) \cdots$ can be $1$.  One has 
\[
\dim R = \rank P + 1.
\]

It follows that if $R$ is an ASL domain on $P$ over $K$, then $P$ possesses a unique minimal element.  In fact, if $\alpha$ and $\beta$ are minimal elements of $P$ with $\alpha \neq \beta$, then the (ASL\,-2) relation of $\varphi(\alpha) \varphi(\beta)$ must be $\varphi(\alpha) \varphi(\beta) = 0$.

We say that an ASL $R$ on $P$ over $K$ is {\em graded} if $R = \bigoplus_{n=0}^{\infty} R_n$ is a noetherian graded algebra over a field $R_0=K$ and $\varphi(P) \subset \bigoplus_{n=1}^{\infty} R_n$.  Furthermore, a graded ASL $R = \bigoplus_{n=0}^{\infty} R_n$ on $P$ over $K$ is called {\em homogeneous} if $\varphi(P) \subset R_1$.

\section{Proof of Theorem \ref{ASLdomain}}

Let $R = \bigoplus_{n=0}^{\infty} R_n$ be a graded ASL on $P$ over $K$ and $S=K[x_{\alpha} : \alpha \in P]$ the polynomial ring in $|P|$ variables over $K$ and define the surjective ring homomorphism $\pi: S \to R$ by setting $\pi(x_\alpha) = \varphi(\alpha)$. The defining ideal $I_R$ of $R=\bigoplus_{n=0}^{\infty} R_n$ is the kernel $\ker(\pi)$ of $\pi$. If $\alpha,\beta \in P$ are incomparable, define
\[
f_{\alpha,\beta}
:= x_\alpha x_\beta - \sum_{i} r_i\,x_{\gamma_{i_1}} x_{\gamma_{i_2}}\cdots , \, \, \, 0 \neq r_i \in K, \, \, \, \gamma_{i_1}\leq \gamma_{i_2}\leq \cdots,
\]
arising from (ASL-2). Then $f_{\alpha,\beta} \in I_R$. Let $\Gc_R$ denote the set of those polynomials $f_{\alpha,\beta}$ for which  $\alpha$ and $\beta$ are incomparable in $P$. Let $<_{\mathrm{rev}}$ denote the reverse lexicographic order \cite[Example~2.1.2(b)]{HHgtm260} on $S$ induced by an ordering of the variables for which $x_{\alpha} <_{\mathrm{rev}} x_{\beta}$ if $\alpha < \beta$ in $P$.  One has each 
\[
{\rm in}_{<_{\mathrm{rev}}}(f_{\alpha, \beta}) = x_\alpha x_\beta.
\]
It follows from (ASL-1) that $\Gc_R$ is a Gröbner basis of $I_R$ with respect to $<_{\mathrm{rev}}$.  In particular, $I_R$ is generated by $\Gc_R$.

\begin{Lemma}
\label{KS}
    Let $R = \bigoplus_{n=0}^{\infty} R_n$ be a graded ASL domain on $P$ over $K$ and $T$ a unique minimal element of $P$.  Then $P\setminus \{T\}$ is pure and strongly connected.
\end{Lemma}

\begin{proof}
    Since the initial ideal ${\rm in}_{<_{\rm rev}}(I_R)$ of $I_R$ with respect to $<_{\rm rev}$ is generated by those monomials $x_\alpha x_\beta$, where $\alpha$ and $\beta$ are incomparable in $P$.  Hence the initial complex \cite[p.~366]{KS} of $I_R$ with respect to $<_{\rm rev}$ is the order complex \cite[p.~119]{HIBIred} of $P$.  Since $I_R$ is a prime ideal, it follows from \cite[Theorem 1]{KS} that $P$ is pure and strongly connected.  Since $T$ is a unique minimal element of $P$, it follows that $P \setminus \{T\}$ is pure and strongly connected.
    \hspace{9.3cm}
\end{proof}

Now, suppose that $R$ is a three-dimensional graded ASL domain on $P$ over $K$ and $T$ is a unique minimal element of $P$.  It follows from Lemma \ref{KS} that $P \setminus \{T\}$ is connected.  Since $\dim R = 3$, one has $\rank (P \setminus \{T\}) = 1$.  Thus $P \setminus \{T\}$ is Cohen--Macaulay (\cite[Exercise (23.2)]{HIBIred}).  In particular, $P$ is Cohen--Macaulay.  Finally, the fundamental theorem on ASL \cite[Theorem (40.1)]{HIBIred} guaranteed that $R$ is Cohen--Macaulay.  This completes a proof of Theorem \ref{ASLdomain}.

\section{Stanley--Reisner rings}
Let $\Delta$ be a simplicial complex on the vertex set $\{x_1, \ldots, x_n\}$ and $S=K[x_1, \ldots, x_n]$ the polynomial ring in $n$ variables over a field $K$.  Let $I_\Delta$ denote the Stanley--Reisner ideal of $\Delta$ and $K[\Delta]:=S/I_\Delta$ the Stanley--Reisner ring.  Let $P_\Delta$ denote the poset consisting of all nonempty faces of $\Delta$, ordered by reverse inclusion.  Since the set of those monomials $u=x_1^{a_1} \cdots x_n^{a_n}$ with $\supp(u) := \{x_i : a_i > 0\} \in \Delta$ is a $K$-basis of $K[\Delta]$, it follows that $K[\Delta]$ is an ASL on $P_\Delta$ over $K$ for which the right-hand side of each of its straightening relations is the empty sum $(=0)$. 

Let $S[t] = K[x_1, \ldots, x_n, t]$ denote the polynomial ring in $n+1$ variables over a field $K$ with each $\deg x_i = e_i > 0$ and $\deg t = e > 0$.  Fix squarefree monomials $w_1, \ldots, w_s$ in $x_1, \ldots, x_n$ together with homogeneous polynomials of the form
\[
f_i = w_i - t g_i(x_1,\ldots, x_n, t), \quad 1 \leq i \leq s 
\] 
for which each monomial appearing in $g_i(x_1,\ldots, x_n, t)$ can be divided by none of the monomials $w_1, \ldots, w_s$, and introduce the quotient ring 
\[
R = S[t]/(f_1, \ldots, f_s).
\]
Let $\Delta$ denote the simplicial complex on $\{x_1, \ldots, x_n\}$ with $K[\Delta] = S/(w_1, \ldots, w_s)$ and $\hat{P}_\Delta = P_\Delta \cup \{T\}$, where $T$ is a unique minimal element of $\hat{P}_\Delta$.   

\begin{Lemma}
\label{gradedASL}
   Let $<$ be a monomial order on $S[t]$ and suppose that $\{f_1, \ldots, f_s\}$ is a Gr\"obner basis of the ideal $(f_1, \ldots, f_s)$ with respect to $<$ with each ${\rm in}_<(f_i) = w_i$.  Then $R$ is a graded ASL on $\hat{P}_\Delta$ over $K$. 
\end{Lemma}

\begin{proof}
Since $\{f_1, \ldots, f_s\}$ is a Gr\"obner basis of $(f_1, \ldots, f_s)$ with respect to $<$ and each ${\rm in}_<(f_i) = w_i$, it follows that the set $\Omega$ of those monomials $u$ of $S[t]$ which is divided by none of $w_1, \ldots, w_s$ is a $K$-basis of $R$.  With the canonical injection $\varphi:\hat{P}_\Delta \hookrightarrow R$ with $\varphi(T) = t$, based on the observation of Stanley--Reisner rings done as above, the condition of (ASL\,-1) that the set of standard monomials is a $K$-basis of $R$ is satisfied.  We show that the condition of (ASL\,-2) is satisfied.  Let $\alpha$ and $\beta$ be nonempty faces of $\Delta$ for which $\alpha$ and $\beta$ are incomparable in $\hat{P}_\Delta$.  
\begin{itemize}
    \item If $\alpha \cup \beta$ is a face of $\Delta$, then
\begin{eqnarray}
\label{aaaaa}
\varphi(\alpha)\varphi(\beta) = \varphi(\alpha \cup \beta) \varphi(\alpha \cap \beta),
\end{eqnarray}
where $\varphi(\alpha \cap \beta)=1$ if $\alpha \cap \beta = \emptyset$.  Since $\varphi(\alpha \cup \beta) < \varphi(\alpha)$ and $\varphi(\alpha \cup \beta) < \varphi(\beta)$ in $\hat{P}_\Delta$, it follows that (\ref{aaaaa}) is the (ASL\,-2) relation of $\varphi(\alpha)\varphi(\beta)$.
\item 
If $\alpha \cup \beta \not\in \Delta$, then $\varphi(\alpha)\varphi(\beta)$ must be divided by one of the monomials $w_1, \ldots, w_s$, say $\varphi(\alpha)\varphi(\beta) = uw_i$, where $u$ is a monomial in $x_1, \ldots, x_n$.  Thus
\begin{eqnarray}
\label{bbbbb}
\varphi(\alpha)\varphi(\beta) = t u g_i(x_1, \ldots, x_n, t).
\end{eqnarray}  
 In the right-hand side of (\ref{bbbbb}), if $v$ is a monomial appearing in $g_i(x_1, \ldots, x_n, t)$ and if $uv$ is divided by one of the monomials $w_1, \ldots, w_s$, say $uv=u'w_j$, where $u'$ is a monomial in $x_1, \ldots, x_n, t$, then 
 \[
 uv=tu'g_j(x_1, \ldots, x_n, t),
 \]
 where 
 \[
 \deg u'g_j(x_1, \ldots, x_n, t)< \deg ug_i(x_1, \ldots, x_n, t).  
 \]
 Continuing with the above procedure, the right-hand side of (\ref{bbbbb}) can finally be expressed as the sum of standard monomials and the final expression is the (ASL\,-2) relation of $\varphi(\alpha)\varphi(\beta)$.
\end{itemize}
This completes a proof of Lemma \ref{gradedASL}. 
\hspace{7cm}
\end{proof}

\section{Proof of Theorem \ref{counterexample}}
We now construct a graded ASL domain of dimension $4$ which is not Cohen--Macaulay over an arbitrary field.  Let $S=K[a,b,c,d,e,t]$ denote the polynomial ring in $6$ variables over a field $K$ with
\[
\deg t = \deg d = \deg e = 1, \quad \deg b = \deg c = 2, \quad \deg a = 3
\]
and $R=S/(f_1,f_2,f_3,f_4)$, where 
\begin{eqnarray*}
f_1&=&bc-t^2de,\\
f_2&=&abe-tb^2d-t^3e^3,\\
f_3&=&acd-tc^2e-t^3d^3,\\
f_4&=&ade-tbd^2-tce^2.
\end{eqnarray*}
Each of $f_1,f_2,f_3,f_4$ is homogeneous with
\[
\deg f_1 = 4, \quad \deg f_2 = 6, \quad \deg f_3 = 6, \quad \deg f_4 = 5.
\]

\begin{Remark}
{\em
    We obtain $f_1,f_2,f_3,f_4$ by weighted homogenization of the defining equations of the non-Cohen--Macaulay domain of \cite[Example 3.14]{CV}.
    }
\end{Remark}
Buchberger's criterion \cite[Theorem 2.3.2]{HHgtm260} guarantees that 
\begin{Lemma}
    \label{initial}
    Let $<_{\rm rev}$ denote the reverse lexicographic order on $S$ induced by the ordering $t < a < b < c < d < e$.  Then $\{f_1,f_2,f_3,f_4\}$ is a Gr\"obner basis of the ideal $(f_1,f_2,f_3,f_4)$ with respect to $<_{\rm rev}$ and
    \[
{\rm in}_{<_{\rm rev}}(f_1)=bc, \quad {\rm in}_{<_{\rm rev}}(f_2)=abe, \quad {\rm in}_{<_{\rm rev}}(f_3)=acd, \quad {\rm in}_{<_{\rm rev}}(f_4)=ade.
    \]
\end{Lemma}

Let $\Delta$ denote the simplicial complex on $\{a,b,c,d,e\}$ with
\[
K[a,b,c,d,e]/(bc,abe,acd,ade) = K[\Delta].
\]

\begin{Corollary}
\label{ASL}
The graded ring $R$ is a graded ASL on ${\hat P}_\Delta$ over $K$. 
\end{Corollary}

The facets of $\Delta$ are 
\[
\{a,b,d\}, \quad \{a,c,e\}, \quad \{b,d,e\}, \quad \{c,d,e\}.
\]
Since ${\rm link}_\Delta(\{a\})$ is disconnected, Reisner's criterion \cite[Theorem (22.3)]{HIBIred} guarantees that $\Delta$ cannot be Cohen--Macaulay over $K$.  In other words, 

\begin{Lemma}
    \label{complex}
The poset ${\hat P}_\Delta$ cannot be Cohen--Macaulay over $K$.  
\end{Lemma}

It then follows from \cite{CV} that

\begin{Corollary}
\label{CM}
The graded ASL $R$ cannot be Cohen--Macaulay.
\end{Corollary}

\begin{Lemma}
    \label{domain}
The graded ASL $R$ is an integral domain.
\end{Lemma}

    \begin{proof}
Let $I=(f_1,f_2,f_3,f_4)$.  Since $R$ is an ASL on ${\hat P}_\Delta$ and $t$ is a unique minimal element of ${\hat P}_\Delta$, it follows from \cite[Corollary 3.6]{Eis} that $t$ is a nonzerodivisor on $R$.

We claim that $b$ is also a nonzerodivisor on $R$. Suppose that
$bh\in I$ for some $h \in S$ such that $h \notin I$. We may assume that $h$ is in its normal form modulo $I$.
The relation $t^3e^3=b(ae-tbd)$ in $R$ gives $t^3e^3h\in I$.
Cancelling $t^3$ in $R$, we obtain $e^3h\in I$. If $h\ne0$,
then $bh,e^3h\in I$ imply that $b \ini (h),e^3 \ini (h)\in J$ and 
\[
  \ini (h)\in (J:b)\cap(J:e^3)
  =(c,ae)\cap(bc,ab,ad)=J,
\]
contrary to the choice of $h$. Thus $b$ is a nonzerodivisor on $R$
and the localization map $R\longrightarrow R[(bt)^{-1}]$ is injective.

Set $A:=K[a,e,b^{\pm1},t^{\pm1}]$. Inverting $bt$ makes both $b$
and $t$ units, so the localization of $S$ at $bt$ is $A[c,d]$.
Since $f_1=bc-t^2de$ and $b$ is a unit in $A$, substitution of
$c=t^2de/b$ gives an isomorphism $A[c,d]/(f_1)\cong A[d]$.
Under this substitution, $f_2$ remains unchanged, while $f_3$ and
$f_4$ become $(t^2d^2/b^2)f_2$ and $(d/b)f_2$, respectively.
Thus the images of $f_2,f_3,f_4$ generate the principal ideal $(f_2)$
in $A[d]$, and hence
$$
R[(bt)^{-1}]\cong A[d]/(f_2).
$$

Now set $d_0:=(abe-t^3e^3)/(tb^2)$. Since $tb^2$ is a unit in $A$,
we have $d_0\in A$ and $f_2=-tb^2(d-d_0)$. Thus $(f_2)=(d-d_0)$
in $A[d]$. Evaluation at $d=d_0$ fixes $A$ and has kernel $(d-d_0)$,
so it induces an isomorphism $A[d]/(f_2)\cong A$. Consequently,
$R[(bt)^{-1}]\cong K[a,e,b^{\pm1},t^{\pm1}]$ is an integral domain.
Since the localization map $R\longrightarrow R[(bt)^{-1}]$ is
injective, $R$ is an integral domain.
\hspace{5.3cm}
\end{proof}

Finally, Lemma \ref{domain} together with Corollaries \ref{ASL} and \ref{CM} proves Theorem \ref{counterexample}.


\begin{thebibliography}{99}
\bibitem{Bac}
K.~Baclawski, Cohen--Macaulay ordered sets, {\em J. of Algebra} {\bf 63} (1980), 226--258.
\bibitem{Bjo}
A.~Bj\"orner, Shellable and Cohen--Macaulay partially ordered sets, {\em Trans. Amer. Math. Soc.} {\bf 260} (1980), 159--183.
\bibitem{CV}
A.~Conca and M.~Varbaro, Square-free Gr\"obner degenerations, {\em Invent. Math.} {\bf 221} (2020), 713--730.
\bibitem{Eis}
D.~Eisenbud, Introduction to algebras with straightening laws, {\em in} ``Ring Theory and Algebra III'' (B.~R.~McDonald, Ed.), Proc. of the third Oklahoma Conf., Lect. Notes in Pure and Appl. Math. No. 55, Dekker, 1980, 243--268.
\bibitem{HHgtm260}
J.~Herzog and T.~Hibi, ``Monomial Ideals'', GTM 260, Springer, 2011. 
\bibitem{HIBIred}
T.~Hibi, ``Algebraic Combinatorics on Convex Polytopes,'' Carslaw Publications, Glebe, N.S.W., Australia, 1992.
\bibitem{HW1}
T.~Hibi and K.-i.~Watanabe, Study of three-dimensional algebras with straightening laws which are Gorenstein domains I, {\em Hiroshima Math. J.} {\bf 15} (1985), 27--54. 
\bibitem{Hoc}
M.~Hochster, Rings of invariants of tori, Cohen--Macaulay rings generated by monomials, and polytopes, {\em Ann. of Math.} {\bf 96} (1972), 318--337.
\bibitem{KS}
M.~Kalkbrener and B.~Sturmfels, Initial complexes of prime ideals, {\em Advances in Math.} {\bf 116} (1995), 365--376.
\bibitem{Rei}
G.~Reisner, Cohen--Macaulay quotients of polynomial rings, {\em Advances in Math.} {\bf 21} (1976), 30--49.
\bibitem{UBC}
R.~Stanley, The upper bound conjecture and Cohen--Macaulay rings, {Stud. Appl. Math.} {\bf 54} (1975), 135--142.
\bibitem{StaGREEN}
R.~P.~Stanley, ``Combinatorics and Commutative Algebra, Second Ed.,'' Birkh\"auser, Boston, 1996.
\bibitem{Terai}
N.~Terai, Some remarks on algebras with straightening laws, {\em J. Pure and Appl. Algebra} {\bf 95} (1994), 87--101.
\end{thebibliography}
\end{document}